\documentclass[3p]{elsarticle}
\usepackage{epsfig}
\usepackage{color}
\usepackage{subfigure}
\usepackage{xcolor}
\usepackage{soul}
\usepackage{amsmath,amssymb}
\usepackage{bm}

\makeatletter
\def\ps@pprintTitle{%
 \let\@oddhead\@empty
 \let\@evenhead\@empty
 \def\@oddfoot{\centerline{\thepage}}%
 \let\@evenfoot\@oddfoot}
 \makeatother

\newcommand{\underbracedmatrix}[2]{%
  \left(\;
  \smash[b]{\underbrace{
    \begin{matrix}#1\end{matrix}
  }_{#2}}
  \;\right)
  \vphantom{\underbrace{\begin{matrix}#1\end{matrix}}_{#2}}
}

\begin{document} 

\title{Travelling waves in a minimal go-or-grow model of cell invasion}

\author[1]{Carles Falc\'o}
\ead{falcoigandia@maths.ox.ac.uk}
\author[1]{Rebecca~M.~Crossley} 
\ead{rebecca.crossley@maths.ox.ac.uk}
\author[1]{Ruth E. Baker}
\ead{ruth.baker@maths.ox.ac.uk}

\cortext[cor1]{Corresponding author.}

\affiliation[1]{organization={Mathematical Institute},
addressline={University of Oxford},
postcode={OX2 6GG},
city={Oxford},
country={United Kingdom}}

\begin{frontmatter}
\begin{abstract} 
We consider a minimal go-or-grow model of cell invasion, whereby cells can either proliferate, following logistic growth, or move, via linear diffusion, and phenotypic switching between these two states is density-dependent. Formal analysis in the fast switching regime shows that the total cell density in the two-population go-or-grow model can be described in terms of a single reaction-diffusion equation with density-dependent diffusion and proliferation. Using the connection to single-population models, we study travelling wave solutions, showing that the wave speed in the go-or-grow model is always bounded by the wave speed corresponding to the well-known Fisher-KPP equation.

\end{abstract} 


\end{frontmatter}

\section{Introduction}

Reaction-diffusion equations are often used to describe the dynamics of invading species, such as cells, and often exhibit travelling wave solutions. 
A classic example of this is the Fisher-Kolmogorov–Petrovsky–Piskunov (FKPP) equation, which assumes that cells move via linear diffusion and proliferate following logistic growth \cite{fisher_wave_1937, kolmogorov1937study}. 
In the last decade however, so-called go-or-grow models \cite{gerlee2012impact,gerlee2016travelling, stepien2018traveling, tursynkozha2023traveling,crossley2024phenotypic} have been suggested to describe the dichotomous behaviour observed in experiments whereby an individual cell either proliferates slowly and migrates quickly, or proliferates quickly and migrates slowly \cite{stein2007mathematical, hoek2008vivo}. 
Under the go-or-grow hypothesis, cells are also able to transition between the migratory and proliferative states via a phenotypic switch. 
The mechanism underlying this switching behaviour is not yet fully understood \cite{tosh2002cells,quesenberry2010cellular}, but could be due to a variety of factors such as surrounding cell density, pressure or environmental features, amongst others. One way to model heterogeneous populations of cells under the go-or-grow hypothesis is to use systems of reaction-diffusion equations.
For specific forms of the phenotypic switching functions, previous works study approximations of the travelling wave solutions by combining the method introduced by Canosa \cite{canosa} with qualitative theory and numerical simulations \cite{stepien2018traveling, tursynkozha2023traveling}. However, these studies also revealed notable differences between the minimum travelling wave speed predicted analytically and the speed observed numerically. 
In addition, minimal progress has been made to understand how invasion phenomena depend on the switching functions more generally \cite{crossley2024phenotypic}.

In this work, we study invasion phenomena in a minimal go-or-grow model, whereby cells either proliferate or migrate, and switching between phenotypic states occurs via general switching functions.
After presenting this model in Section~\ref{sec:model}, we first demonstrate, via a formal fast reaction limit, the connection between a model involving heterogeneous populations and a model consisting of a single population of cells, with density-dependent diffusion and proliferation terms. 
Interestingly, this reveals a connection between multiple phenotype models with density-dependent switching functions, and degenerate diffusion models consisting of a single phenotype, for which the standard linear stability analysis does not yield accurate travelling wave speed predictions \cite{murray2001mathematical, el2021travelling}.
We then proceed in Section~\ref{sec:TWAhetero} to perform a travelling wave analysis of the minimal go-or-grow model presented in this work and show that the minimum travelling wave speed emitted from the model does not exceed that of the FKPP equation. 
Furthermore, we demonstrate that for specific forms of the phenotypic switching function, solutions of the multiple phenotype model do not attain this minimum travelling wave speed. 
This provides an analytical explanation for the observed discrepancies  in prior studies between the minimum wave speed predicted analytically and the numerically observed travelling wave speed \cite{stein2007mathematical}.

\section{A minimal go-or-grow model}\label{sec:model}

We consider a heterogeneous cell population consisting of cells which are either proliferative or invasive. The density of the proliferative population is represented by $\rho_1(\mathbf{x},t)$ and the density of the invasive population is given by $\rho_2(\mathbf{x},t)$, for $(\mathbf{x},t)\in\mathbb{R}^d\times[0,+\infty)$. 
For the sake of simplicity, we assume that invasive cells move via linear diffusion, and can switch to the proliferative phenotype at a rate given by the switching function $\Gamma_1(\rho)$, which depends on the total cell density $\rho = \rho_1+\rho_2$. 
On the other hand, we assume that the proliferative cells follow logistic growth, and can switch to the invasive phenotype at a rate given by $\Gamma_2(\rho)$. 
With this, the model reads
\begin{equation}\label{eq: full model}
\begin{split}
    \partial_t\rho_1 & = \Delta\rho_1 - \rho_1\Gamma_1(\rho) + \rho_2\Gamma_2(\rho)\,,
    \\
    \partial_t\rho_2 & = \rho_2\left(1-{\rho}\right) + \rho_1\Gamma_1(\rho) - \rho_2\Gamma_2(\rho)\,.
\end{split}
\end{equation}
Following from standard rescaling of space and time, the diffusion coefficient and the intrinsic proliferation rate in Eqs.~\eqref{eq: full model} have been set to one without loss of generality. 
We also note that this model does not account for cell death, however apoptosis can be readily incorporated into the model, as demonstrated in \cite{stepien2018traveling,tursynkozha2023traveling}, without altering the analysis. 
The results in this work generalise directly to this scenario.

The model given in Eqs.~\eqref{eq: full model} is studied in \cite{stepien2018traveling, tursynkozha2023traveling} for specific choices of switching functions, which facilitate analytic tractability. 
Here, and unless stated otherwise, we keep $\Gamma_1(\rho)$ and $\Gamma_2(\rho)$ general. 
Following the biological intuition that regions of high cell density lead to contact inhibition or nutrient limitation, and hence to inhibition of proliferation \cite{puliafito2012collective, falco2024mechanical}, typical switching functions are based on the assumption that $\Gamma_1(\rho)$ is non-increasing, and $\Gamma_2(\rho)$ is non-decreasing. 
While the model presented here is simpler, it remains closely related to those examined in \cite{gerlee2012impact,gerlee2016travelling}, which consider density-dependent diffusion, and in \cite{crossley2024phenotypic}, which further incorporates the influence of the extracellular matrix on cell invasion. 

\subsection{Connecting heterogeneous to single-population models}

We investigate the fast phenotypic switching limit of Eqs.~\eqref{eq: full model}. 
In other words, we assume $\Gamma_1(\rho)= \Tilde{\Gamma}_1(\rho)/\epsilon$, and $\Gamma_2(\rho)= \Tilde{\Gamma}_2(\rho)/\epsilon$ where $\Tilde{\Gamma}_1(\rho), \, \Tilde{\Gamma}_2(\rho) \sim O(1)$ for $\epsilon\ll 1$.
Upon formally taking the limit $\epsilon\rightarrow 0$ in Eqs.~\eqref{eq: full model}, we obtain $\rho_1\Tilde{\Gamma}_1(\rho) = \rho_2\Tilde{\Gamma}_2(\rho)$. 
Consider now the equation governing the evolution of the total cell density over time,
\begin{equation}
    \partial_t \rho = \Delta\rho_1 + \rho_2(1 - \rho)\,.\label{eq: total density pde}
\end{equation}
In the limit $\epsilon\rightarrow 0$, the two cell densities $\rho_1$ and $\rho_2$ can be related to the total cell density via
\begin{equation}
    \rho_1 = \frac{\Gamma_2(\rho)}{\Gamma_1(\rho) + \Gamma_2(\rho)}\,\rho\,,\quad\rho_2 = \frac{\Gamma_1(\rho)}{\Gamma_1(\rho) + \Gamma_2(\rho)}\,\rho\,.\label{eq: fractions}
\end{equation}
Combining Eqs.~\eqref{eq: total density pde} and~\eqref{eq: fractions} we obtain the following equation for the evolution of the total cell density
\begin{equation}
    \partial_t\rho = \nabla\cdot\left(D(\rho)\nabla\rho\right) +  r(\rho)\rho\,, \label{eq: total density limit}
\end{equation}
with the density-dependent diffusion coefficient $D(\rho)$, and proliferation rate $r(\rho)$ given by
\begin{equation*}
    D(\rho) = \frac{\Gamma_1(\rho)\Gamma_2'(\rho)-\Gamma_1'(\rho)\Gamma_2(\rho)}{\left(\Gamma_1(\rho) + \Gamma_2(\rho)\right)^2}\,\rho + \frac{\Gamma_2(\rho)}{\Gamma_1(\rho) + \Gamma_2(\rho)},\qquad r(\rho) = \frac{\Gamma_1(\rho)(1-\rho)}{\Gamma_1(\rho) + \Gamma_2(\rho)}\,.
\end{equation*}
In this case, assuming that $\Gamma_1$ is a non-increasing function of the total cell density is sufficient to show that the diffusion coefficient does not become negative, i.e. $D(\rho)\geq 0$ for any $\rho\geq 0$; and also that the per capita proliferation rate is a non-increasing function of density, i.e. $r'(\rho)\leq 0$. 
As mentioned earlier, this assumption on $\Gamma_1$ aligns with basic biological intuition. 
Moreover, under the assumption that cells remain proliferative at low densities, so that $\Gamma_2(0) = 0$, we obtain that the diffusion coefficient follows $D(\rho)\sim 2\Gamma_2'(0)/\Gamma_1(0)\rho + O(\rho^2)$, which resembles typical non-linear equations with degenerate diffusion.

Finally, we highlight that in the case of constant switching rates ($\Gamma_1'=\Gamma_2'=0$) we recover the well-known FKPP equation \cite{murray2001mathematical}
\begin{equation}
    \partial_t\rho = \theta \Delta\rho + (1-\theta)\rho(1-\rho)\,,\label{eq: FKPP}
\end{equation}
where $\theta\in(0,1)$ is determined by the constant ratio $\Gamma_1/\Gamma_2$. 
Equations of this type often exhibit travelling wave solutions $\rho(\mathbf{x},t) = U(\mathbf{r}\cdot\mathbf{x}-ct)$, where $\mathbf{r}$ is the direction of propagation, and $c$ corresponds to the wave speed. 
In particular, the FKPP equation (Eq.~\eqref{eq: FKPP}) admits travelling wave solutions with constant speed $c\geq c_{\text{min}} = 2\sqrt{\theta(1-\theta)}$. 
We note that $c_{\text{min}}\leq 1$ for $\theta\in(0,1)$, which we use as a reference to compare with the speed of propagation of travelling wave solutions to Eqs.~\eqref{eq: full model}.

\section{Travelling wave analysis for heterogeneous populations} \label{sec:TWAhetero}
\begin{figure}[b]
    \centering
    \includegraphics[width = \textwidth]{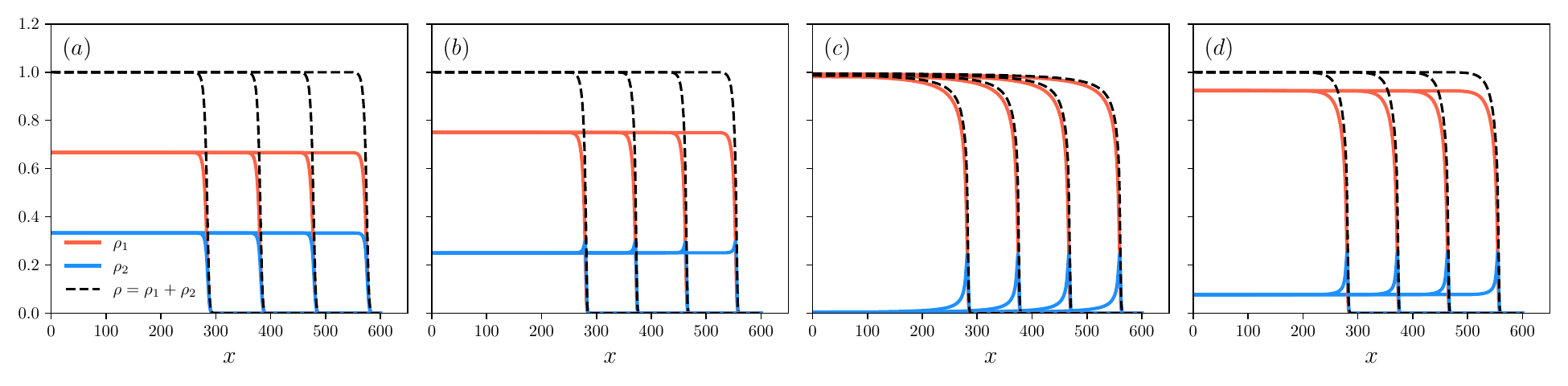}
    \caption{Numerical solutions to Eqs.~\eqref{eq: full model} subject to the initial conditions: $\rho_1(x,0) =\rho_2(x,0) = 0.2$ for $x<100$, and $\rho_1(x,0) = \rho_2(x,0) = 0$ for $x\geq 100$, and for different choices of the switching functions $\Gamma_1,\,\Gamma_2$. Solutions are plotted at times $t = 200, 300, 400, 500$. Chosen switching functions are as follows: (a) $\Gamma_1(\rho) = 0.5,\,\Gamma_2(\rho) = 1$; (b) $\Gamma_1(\rho) = 0.5,\,\Gamma_2(\rho) = 1.5\rho$; (c) $\Gamma_1(\rho) = 0.5(1-\rho),\,\Gamma_2(\rho) = 1.5\rho$; (d)$\Gamma_1(\rho) = 0.5(1-\rho^2/(0.5^2+\rho^2)),\,\Gamma_2(\rho) = 1.5\rho^2/(0.5^2+\rho^2)$.}
    \label{fig: tws}
\end{figure}

We now revisit Eqs.~\eqref{eq: full model} to analyse travelling wave solutions of the form: $\rho_1(\mathbf{x},t) = U_1(\mathbf{r}\cdot\mathbf{x}-ct)$, $\rho_2(\mathbf{x},t) = U_2(\mathbf{r}\cdot\mathbf{x}-ct)$.
Fig.~\ref{fig: tws} shows travelling wave solutions of Eqs.~\eqref{eq: full model} subject to different choices of the switching functions. 
By defining $z = \mathbf{r}\cdot\mathbf{x}-ct\in(-\infty,+\infty)$, and writing Eqs.~\eqref{eq: full model} in travelling wave coordinates we obtain
\begin{equation}\label{eq: tw coordinates}
\begin{split}
     c\frac{\mathrm{d}U_1}{\mathrm{d}z} + \frac{\mathrm{d}^2U_1}{\mathrm{d}z^2} - U_1\Gamma_1(U)+U_2\Gamma_2(U)&=0\,,
    \\ 
     c\frac{\mathrm{d}U_2}{\mathrm{d}z} + U_2(1-U) + U_1\Gamma_1(U)-U_2\Gamma_2(U)&=0\,,
\end{split}
\end{equation}
where $U(z) = U_1(z) + U_2(z)$. 
Next, we look at steady states of the system~\eqref{eq: tw coordinates}.

First note that the origin $(U_1,U_2) = (0,0)$ is always a steady state, and we assume $U_1,U_2\rightarrow 0$ as $z\rightarrow+\infty$. 
The other steady states, corresponding to $z\rightarrow-\infty$, depend on the choice of phenotypic switching functions (see Fig.~\ref{fig: tws}) and belong to the set $\left\{(U_1,U_2)\in[0,1]^2:U_1\Gamma_1(U)=U_2\Gamma_2(U),\;U_2(1-U) = 0\right\}$. 
In particular, a steady state with $U_2=0$ can only occur if there exists a $U^*>0$ such that $\Gamma_1(U^*) = 0$. 
In this case, we obtain
\begin{equation*}
    (U_1,U_2) = (U^*,0)\,.
\end{equation*}
If $\Gamma_1(U)>0$ for all $U>0$, then we obtain
\begin{equation*}
    (U_1,U_2) = \left(\frac{\Gamma_2(1)}{\Gamma_1(1)+\Gamma_2(1)}, \frac{\Gamma_1(1)}{\Gamma_1(1)+\Gamma_2(1)}\right)\,.
\end{equation*}

\subsection{Dispersion relation and minimum travelling wave speed}
\begin{figure}
    \centering
    \includegraphics[width = .9\textwidth]{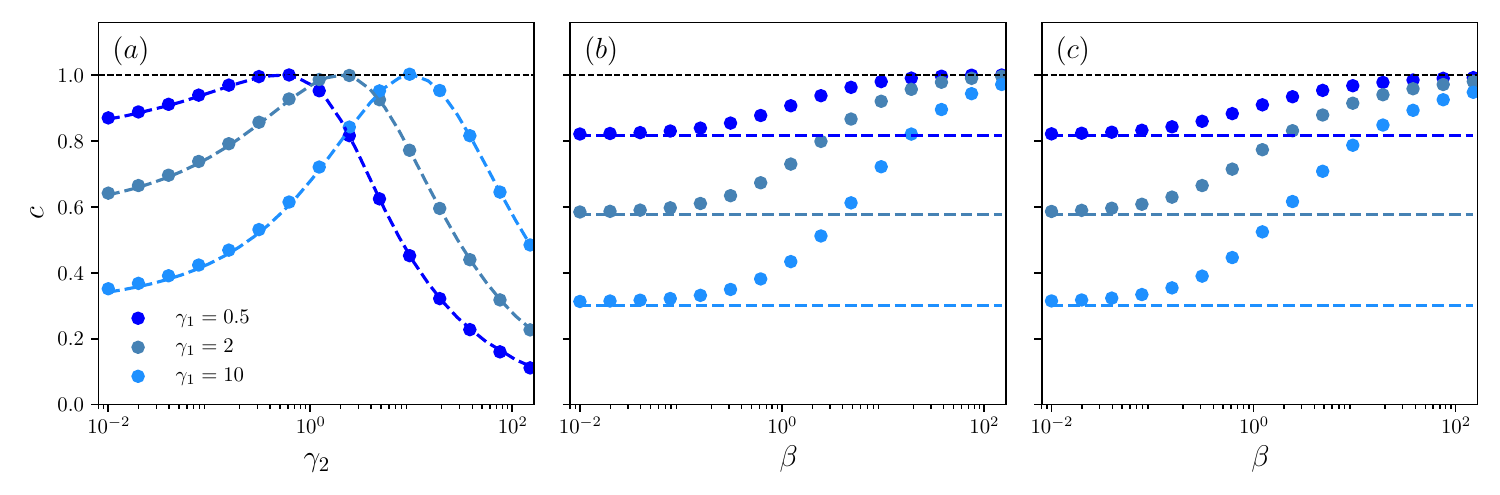}
    \caption{Relationship between the estimated numerical speed (dots) of travelling wave solutions to Eqs.~\eqref{eq: full model}, and the minimum wave speed predicted by Eq.~\eqref{eq: csigma} (dashed lines). Numerical wave speed estimated using a trapezoidal rule applied to the expression $c = \int_{-\infty}^{+\infty} U_2(z)(1-U(z))\,\mathrm{d}z$, which follows from Eq.~\eqref{eq: tw coordinates}. Eqs.~\eqref{eq: full model} are solved on a domain of length 7000, for a maximum simulation time of 6500, and initial conditions given by $\rho_1(x,0)=\rho_2(x,0) = 0.2$ if $x<100$ and $\rho_1(x,0)=\rho_2(x,0) = 0$ if $x\geq 100$. Chosen switching functions are as follows: (a) $\Gamma_1(\rho) = \gamma_1,\,\Gamma_2(\rho) = \gamma_2$; (b) $\Gamma_1(\rho) = \gamma_1,\,\Gamma_2(\rho) = \beta\rho$; (c) $\Gamma_1(\rho) = \gamma_1(1-\rho^2/(0.5^2+\rho^2)),\,\Gamma_2(\rho) = \beta\rho^2/(0.5^2+\rho^2)$.}
    \label{fig: speeds}
\end{figure}

When seeking an expression for the travelling wave speed, $c$, it is often commonplace to perform a linear stability analysis centred around the origin, which often yields a lower bound on the speed: $c\geq c_{\text{min}}$, for some minimum wave speed $c_{\text{min}}\geq 0$.
Unfortunately, for the model presented in this work, analytical progress following this approach is limited and instead we shift the investigation to identify the dispersion relation, which relates the travelling wave speed with the dynamics at the leading front. 
We observe that at the leading front, cell densities are small, and hence we can linearise Eqs.~\eqref{eq: tw coordinates} to give
\begin{equation} \label{eq: linear tw}
\begin{split}
     c\frac{\mathrm{d}U_1}{\mathrm{d}z} + \frac{\mathrm{d}^2U_1}{\mathrm{d}z^2} - \gamma_1U_1+\gamma_2U_2&=0\,,
    \\ 
     c\frac{\mathrm{d}U_2}{\mathrm{d}z}  + \gamma_1U_1+(1-\gamma_2)U_2&=0\,,
\end{split}
\end{equation}
where we denote $\gamma_i = \Gamma_i(0)$ for $i=1,2$.
By substituting the ansatz $U_1(z)\sim Ae^{-\sigma z}$, $U_2(z)\sim B e^{-\sigma z}$, for $z\gg 1$ and $A, B\geq 0$ constants, we obtain
\begin{equation*}
\underbracedmatrix{
    -c\sigma+\sigma^2-\gamma_1 & \gamma_2
    \\
    \gamma_1 & -c\sigma + 1-\gamma_2}{M}
\begin{pmatrix}
    A \\ B
\end{pmatrix} = 0\,.
\end{equation*}
In order to ensure $A,\,B\neq 0$, we require $\det M = 0$, which gives
\begin{equation*}
    c^2-c\left(\sigma + \frac{1-\gamma_1-\gamma_2}{\sigma}\right)+1-\gamma_2-\frac{\gamma_1}{\sigma^2} = 0\,.
\end{equation*}
The travelling wave speed, $c$, can then be solved for as a function of $\sigma$, giving the dispersion relation:
\begin{equation}\label{eq: csigma}
    c(\sigma) = \frac{1}{2}\left(f(\sigma)+\sqrt{f(\sigma)^2-4(1-\gamma_2)+\frac{4\gamma_1}{\sigma^2}}\,\right), \quad f(\sigma) = \sigma + \frac{1-\gamma_1-\gamma_2}{\sigma}.
\end{equation}
The expression above predicts the existence of a minimum travelling wave speed, $c_{\text{min}}$ for $\sigma = \sigma^*$, such that $c_{\text{min}}=c(\sigma^*)$. 
Finding the minimum speed analytically seems challenging; instead we observe that $c(1) = 1$, and by differentiating $c(\sigma)$ with respect to $\sigma$ and evaluating at $\sigma = 1$, we obtain
\begin{equation*}
    \left(\frac{\mathrm{d}c}{\mathrm{d}\sigma}\right)_{\sigma = 1}= \frac{\gamma_2-\gamma_1}{\gamma_1+\gamma_2}\,.
\end{equation*}
Thus, for $\gamma_1 = \gamma_2$, $\sigma^* = 1$ yields the minimum wave speed, and in this case it is given by $c_\text{min} = 1$. 
On the other hand, when $\gamma_1>\gamma_2$ we have $({\mathrm{d}c}/{\mathrm{d}\sigma})_{\sigma = 1}<0$, and $\sigma^*>1$ with $c_{\text{min}}<1$. 
An analogous argument yields $\sigma^*<1$ and $c_{\text{min}}<1$ for $\gamma_1<\gamma_2$. 
Hence, the minimum travelling wave speed displayed by Eqs.~\eqref{eq: full model} is, at most, the minimum travelling wave speed corresponding to the FKPP model, Eq.~\eqref{eq: FKPP}. Fig.~\ref{fig: speeds}a shows an excellent agreement between the numerically observed speed and the prediction from Eq.~\eqref{eq: csigma} for a fully linear model with constant switching functions. Next, we explore whether this agreement holds for more general nonlinear models.

\subsection{Explicit expressions for the minimum travelling wave speed}

Fully analytical expressions for $c_{\text{min}}$ can be obtained under the assumption that $\gamma_1 = 0$ or $\gamma_2 = 0$. We present first the $\gamma_2 =0$ case as it represents a more biologically realistic scenario where cells remain proliferative at low densities. In this case, Eq.~\eqref{eq: csigma} simplifies to
\begin{equation*}
    c(\sigma) = \frac{1}{2}\left(\sigma + \frac{1-\gamma_1}{\sigma}+\frac{\left|\sigma^2-1-\gamma_1\right|}{\sigma}\right)\,.
\end{equation*}
Note then that $c(\sigma)$ increases for $\sigma > \sqrt{1+\gamma_1}$, and decreases for $\sigma < \sqrt{1+\gamma_1}$. Hence, $\sigma^* = \sqrt{1+\gamma_1}$, and 
\begin{equation}\label{eq: cmin gamma2=0}
    c_{\text{min}} = \frac{1}{\sqrt{1+\gamma_1}}\,.
\end{equation}

\noindent When $\gamma_1 = 0$, a similar argument yields the minimum wave speed
\begin{align*}
    c_{\text{min}} = \begin{cases}
        \sqrt{1-\gamma_2} &\mbox{for } \gamma_2 < 1\,,
        \\
         0  &\mbox{for } \gamma_2 > 1\,.
    \end{cases}
\end{align*}

\subsection{The minimum wave speed is not attainable when $\Gamma_2(0)=0$}
Based on the behaviour observed in Fig.~\ref{fig: speeds}b-c, we focus on the case $\gamma_2=\Gamma_2(0)=0$, with $\Gamma_2\neq0$ to ensure that travelling wave solutions exist.
By linearising Eqs.~\eqref{eq: tw coordinates} subject to $\Gamma_2(0)=0$, representing the region near the invading front where cell densities are low, we obtain
\begin{equation*} 
\begin{split}
     c\frac{\mathrm{d}U_1}{\mathrm{d}z} + \frac{\mathrm{d}^2U_1}{\mathrm{d}z^2} - \gamma_1U_1&=0\,,
    \\ 
     c\frac{\mathrm{d}U_2}{\mathrm{d}z}  + \gamma_1U_1+U_2&=0\,.
\end{split}
\end{equation*}
Solving this system of equations, subject to $U_1, \, U_2 \rightarrow 0$ as $z\rightarrow \infty$, yields
\begin{equation*}
    U_1(z) = Ae^{-\lambda(c)z}, \quad   U_2(z) = \frac{\gamma_1}{c\lambda(c) - 1}\,Ae^{-\lambda(c)z},\quad \lambda(c) = \frac{1}{2}\left(c+\sqrt{c^2+4\gamma_1}\right)\,,
\end{equation*}
for an arbitrary constant $A>0$. Positivity of travelling wave solutions (i.e. $U_1, U_2\geq0$ for all $z\in\mathbb{R}$)
requires $c\lambda(c) - 1>0$, which is only satisfied when $c>c_{\text{min}}$ (as in Eq.~\eqref{eq: cmin gamma2=0}).   
In particular, we obtain that travelling wave solutions for $\Gamma_2(0)=0$ exhibit a wave speed exceeding $c_{\text{min}}$, i.e. $c>c_{\text{min}}$ when $\Gamma_2(0)=0$.

\section{Discussion and open problems}

In this work, we examined a minimal go-or-grow model, where cells can either proliferate following logistic growth or move through linear diffusion. Via a formal fast reaction limit in the fast phenotypic switching regime, we revealed a connection between heterogeneous models and single-population models. This relationship implies an approximate equivalence between discrete-state representations of cell populations and models where diffusion and proliferation rates are determined by density. Moreover, our analysis suggests that the invasion speed in travelling wave solutions of heterogeneous models is bounded by the wave speed of solutions to single-populations models. We highlight that this phenomenon is likely to change when energy or resource constraints associated to diffusion and proliferation are taken into account \cite{crossley2024phenotypic}. More research is needed, however, to better understand how to model such energy constraints and their impact on cell invasion.

Furthermore, our travelling wave analysis also reveals a lower bound on the invasion speed of travelling wave solutions to go-or-grow models. This estimate, based on a linear stability analysis,  breaks down for nonlinear models (see Fig.~\ref{fig: speeds}b-c), as observed in similar models \cite{stepien2018traveling, el2021travelling, colson2021travelling, crossley2023}. Numerically, we observe that $c_{\text{min}}\leq c\leq 1$ for the model given by Eqs.~\eqref{eq: tw coordinates}, and that $c\rightarrow 1$ as the nonlinear switching rates increase. Although this behaviour intuitively makes sense when considering the connection to single-population models, we do not currently have a rigorous proof. We leave these conjectures for further investigation.

\section*{Acknowledgments}
RMC is supported by the Engineering and Physical Sciences Research Council (EP/T517811/1) and the Oxford-Wolfson-Marriott scholarship at Wolfson College, University of Oxford.
This work was supported by a grant from the Simons
Foundation (MP-SIP-00001828, REB).

\end{document}